\documentclass{article}

\usepackage{amssymb}      
\usepackage{amsmath}      
\usepackage{amsthm}      
\usepackage{amsbsy}      
\usepackage{amscd}      
\usepackage[dvips]{graphics}
\usepackage{mathrsfs}
\usepackage{pstricks}

\newtheorem{teo}{Theorem}
\newtheorem{lema}[teo]{Lemma}
\newtheorem{prop}[teo]{Proposition}
\newtheorem{corolar}[teo]{Corollary}

\theoremstyle{definition}
\newtheorem{dfn}[teo]{Definition}
\newtheorem{exemplu}[teo]{Example}

\theoremstyle{remark}
\newtheorem{rem}[teo]{Remark}
\newtheorem{conj}[teo]{Conjecture}

\newcommand{\eq}{{\ =\ }}

\newcommand\er{{\mathbb R}}
\newcommand\ce{{\mathbb C}}
\newcommand\es{{\mathbb S}}
\newcommand\de{{\mathbb D}}

\newcommand\q{{\bf Q}}
\newcommand\f{{\bf F}}
\newcommand\p{{\bf P}}

\newcommand\braid{{\rm B}}

\newcommand\ar{{\rm R}}

\begin{document}
\date{\today}
\title{
Exchange Moves and Fiedler Polynomial\footnote {
This is part of my Ph.D. Thesis I defended at Columbia University in the spring of 2001}
\author{Radu Popescu\\      
\small \em Institute of Mathematics of the Romanian Academy,           
 P.O. Box 1-764, RO-014700, Bucharest, Romania \\ }}
\date{{}}

\maketitle

A link is an embedding of $\amalg_{i=1}^n\es^1$ in $\es^3$ or in $\er^3$. 
A famous problem is the classification  of links up to isotopy. 
This problem has an algebraic solution through braids. A braid is an embedding of 
$\amalg_{i=1}^n I$ into the cylinder $\de^2\times I$, 
such that the intervals are considered to strictly decrease from
the top to the bottom disc. The end points are considered to be fixed. 

The braids form a group, the multiplication being
given by putting the cylinders in top of one another, such that the end points correspond. The identity element is 
the braid formed by all the strings going straight from the point $i\times 1$ to $i\times 0$. 
The braid group $\braid_n$ is finitely presented given the number of strings also called the braid index.
The standard set of generators is $\sigma_1,\dots,\sigma_{n-1} $, $\sigma_i$ being represented in the picture below.
Moreover the braid $\sigma_i^{-1}$ is the one in which the crossing of the $i$-th string
is in front of the the $i+1$-th.

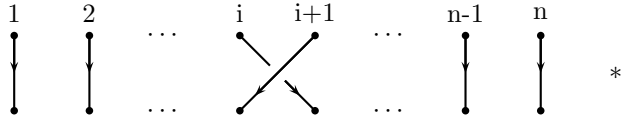
\begin{figure}[ht]
\begin{pspicture}(0,1.7)
\pscircle*(5,1){.05}
\pscircle*(5,0){.05}
\pscircle*(6,1){.05}
\pscircle*(6,0){.05}
\rput(7,1){$\cdots$}
\rput(7,0){$\cdots$}
\pscircle*(8,1){.05}
\pscircle*(8,0){.05}
\pscircle*(9,1){.05}
\pscircle*(9,0){.05}
\rput(10,1){$\cdots$}
\rput(10,0){$\cdots$}
\pscircle*(11,1){.05}
\pscircle*(11,0){.05}
\pscircle*(12,1){.05}
\pscircle*(12,0){.05}
\psline[linewidth=.8pt](5,1)(5,0)
\psline[linewidth=.8pt]{->}(5,1)(5,0.5)
\psline[linewidth=.8pt](6,1)(6,0)
\psline[linewidth=.8pt]{->}(6,1)(6,0.5)
\psline[linewidth=.8pt](9,1)(8,0)
\psline[linewidth=.8pt]{->}(9,1)(8.2,0.2)
\psline[linewidth=.8pt](8,1)(8.4,0.6)
\psline[linewidth=.8pt](8.6,0.4)(9,0)
\psline[linewidth=.8pt]{->}(8.6,0.4)(8.8,0.2)
\psline[linewidth=.8pt](11,1)(11,0)
\psline[linewidth=.8pt]{->}(11,1)(11,0.5)
\psline[linewidth=.8pt](12,1)(12,0)
\psline[linewidth=.8pt]{->}(12,1)(12,0.5)
\rput(5,1.3){1}
\rput(6,1.3){2}
\rput(8,1.3){i}
\rput(9,1.3){i+1}
\rput(11,1.3){n-1}
\rput(12,1.3){n}
\rput(13,.5){$\ast$}
\end{pspicture}
\caption{Generator $\sigma_i$}
\label{generator:sigma_i}
\end{figure}
\vskip .1in

Adding a trivial string one obtains a braid in $\braid_{n+1}$. This give us an embedding 
$\braid_n \hookrightarrow \braid_{n+1}$. 

A braid can be closed to a link connecting every point $i\times 0$ to
$i\times 1$ by an arc which go around an axis perpendicular to the cross-section of the cylinder. The axis is shown
in the above picture by the star. A link obtained in this manner is also called a closed braid.

A theorem of Alexander says that any link in $\es^3$ (or $\er^3$) can be isotoped into a closed braid.
Markov theorem says that the closures of two braids $\alpha,\beta$ are isotopic if and only if the braids are related by a
finite sequence of the following two moves and their inverses. 

\[
\begin{array}{l}

(1)\ \beta\eq \gamma\cdot\alpha\cdot\gamma^{-1},\hskip 1in\alpha,\beta,\gamma\in\braid_n\\
 
(2)\ \beta\eq \alpha\cdot\sigma_n^{\pm 1},  
{\rm\ changing\ the\ braid\ index}, n\leftrightharpoons n+1
\end{array}
\]

The trouble in the classification problem mentioned above is made by the mixture
of the two moves. To be more precise for each isotopy class of links there are infinitely many braids
which close to a link in such a class. But there is a braid representative of minimum braid index representing the link.
The problem is that given a braid with a larger that minimum braid index,
representing a link ( up to isotopy), it might be the case that in order to simplify it ( meaning to decrease
the braid index) one should go up in the tower of braid groups, and then go to a lower index than the one at the beginning.

Examples of this kind were given by Rudolph in ~\cite{Rudolph} for the unlink with two components, refined by Morton in ~\cite{MR84m:57006} to an example for the unknot. 

To overcome the difficulty Birman and Menasco in the paper 
~\cite{MR92g:57010a} described a new move which doesn't change the 
braid index and may be accomplished only by crossing the braid 
axis as in Markov move (2), above.  
This new move is given in the following:

\begin{dfn}
\label{def:exchange}
Two braids $\beta_1, \beta_2$
are said to be related by an \underline {exchange move}, or 
called \underline {exchange related},
if $\beta_1\eq X\sigma_n^{-1} Y\sigma_n\ and\ \beta_2\eq X\sigma_n 
Y\sigma_n^{-1}$, 
where 
$X$ and $Y\in\braid_n$.
\end{dfn}

The move changes the conjugacy class of generic braids for $n\geq 4$. 
Our goal is to discover a simple and definitive test for verifying 
whether two braids which are
related by an exchange move are conjugate or not.
In ~\cite{prep1} the authors describe an algorithmic
solution for the conjugacy problem in $\braid_n$.
This gives, of course, a complete answer to our question,
however the algorithm is complicated and one might hope for a simpler 
solution in our special situation. The matter turns out, however, 
to be subtle. All the easy invariants take an identical value
on exchange-equivalent conjugacy classes.

Let me describe briefly the content of the rest of the paper.
In Section {\bf\ref{section:fiedler}} I give the
definition of the Fiedler invariant, which is a class
invariant. I state its properties in Theorem {\bf\ref{teo:fiedler}}
and give a sketch of the proof for the case of braided knots.
Although implicit in Morton's work ~\cite{MR2000i:57025}, I prove
in Proposition {\bf\ref{prop:fiedlerme}} that Fiedler's invariant
is a type one invariant. Then in Section {\bf\ref{section:excfiedler}} 
I use Fiedler's polynomial to detect exchange related braids. Necessary and sufficient conditions
for this are given in Proposition {\bf\ref{prop:exrelength}}. 
At the end of the section I give explicit examples (see Examples
{\bf\ref{ex:braid1}}---{\bf\ref{ex:braid3}}).

\section{Fiedler's invariant}
\label{section:fiedler}

Fiedler in ~\cite{MR94c:57006} describes in a general setting an isotopy
invariant for knots. He considers knots K embedded in an orientable 
3-manifold E which is the total space of a real line bundle over a surface. 
I will consider in what follows only braided knots. 
The knots are embedded in E$\eq\er^3\setminus\{z-axis\}$ and 
the surface is S$\eq{\er^2}\setminus\{0\}$.

For this particular case we get in fact an invariant 
of closed braids up to conjugacy, as described in Fiedler's paper.
Applying this invariant to $\widehat{\beta_1}$ and 
$\widehat{\beta_2}$, exchange related braided knots,
we see that there are a lot of cases in which it can distinguish 
when they are not conjugate. 

In the sequel I will consider only braids $\gamma = \prod_{r=1}^k{\sigma_{i_r}^{\epsilon_r}}$ 
such that their closures $\hat\gamma$ are knots. By smoothing I will understand the replacement shown 
in the picture below. 

\begin{figure}[ht]
\begin{pspicture}(0,1.7)
\pscircle*(4,1){.05}
\pscircle*(4,0){.05}
\pscircle*(5,1){.05}
\pscircle*(5,0){.05}
\pscircle*(6,1){.05}
\pscircle*(6,0){.05}
\pscircle*(7,1){.05}
\pscircle*(7,0){.05}
\pscircle*(10,1){.05}
\pscircle*(10,0){.05}
\pscircle*(11,1){.05}
\pscircle*(11,0){.05}
\psline[linewidth=.8pt](4,1)(4.4,0.6)
\psline[linewidth=.8pt]{->}(4.6,0.4)(5,0)
\psline[linewidth=.8pt]{->}(5,1)(4,0)
\psline[linewidth=.8pt]{->}(6,1)(7,0)
\psline[linewidth=.8pt](7,1)(6.6,0.6)
\psline[linewidth=.8pt]{->}(6.4,0.4)(6,0)
\psline[linewidth=.8pt]{->}(8,0.5)(9,0.5)
\pscurve[linewidth=.8pt]{->}(10,1)(10.4,0.5)(10,0)
\pscurve[linewidth=.8pt]{->}(11,1)(10.6,0.5)(11,0)
\rput(4,1.3){$i$}
\rput(5,1.3){$i+1$}
\rput(6,1.3){$i$}
\rput(7,1.3){$i+1$}
\rput(10,1.3){$i$}
\rput(11,1.3){$i+1$}
\rput(5.5,0.5){or}
\rput(4.5,-0.5){$\epsilon\eq 1$}
\rput(6.5,-0.5){$\epsilon\eq -1$}
\end{pspicture}
\vskip 0.3in
\caption{Smoothing}
\label{smoothing}
\end{figure}
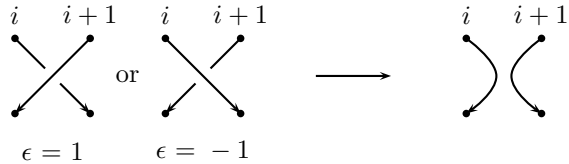
\vskip .1in

\begin{dfn}
\label{dfn:ascending}
The ascending string of a crossing in a given closed braid $\hat\gamma$
(counting all the braid strings from 1 to n at the top of each crossing) is:
\begin{itemize}
\item if the sign of the crossing $\sigma_{i_r}$ is $\epsilon_r = 1$ then the \underline{ascending string} 
is the $i_r$-th string after smoothing the crossing;
\item if the sign of the crossing $\sigma_{i_r}$ is $\epsilon_r = -1$ then the \underline{ascending string} 
is the $i_{r+1}$-th string after smoothing the crossing.
\end{itemize}
\end{dfn}

\begin{dfn}
\label{dfn:fiedler} 
The \underline{Fiedler polynomial} for any 
$\gamma\in\braid_n$ is given by:
\[
\f_{\hat\gamma}(x) =\sum_{r=1}^k\epsilon_r{x^{2m(r)-n}}
\]
where m(r) is the winding number of the ascending string 
around the braid axis, after smoothing the r-th 
crossing of the braid.
\end{dfn}

Let $\pi:\braid_n\longrightarrow\Sigma_n$ be the homomorphism 
which assigns to each braid its associated  permutation. 
Let $s_i$ be the transposition $(i,i+1)$. 
I will denote with $\pi_X$ the image of the braid 
$X$ through $\pi$, so, $\pi(\sigma_i)\eq \pi_{\sigma_i}\eq s_i$.

\vskip 0.7in
\begin{teo}[Fiedler] 
\label{teo:fiedler}
$\f_{\hat\gamma}(x)$ is a conjugacy class invariant of $\gamma \in\braid_n$ 
where $\hat\gamma$ is a knot, and has the following properties:
\begin{enumerate}
\item $\f_{\hat\gamma}(x)$ is a symmetric Laurent polynomial in x.
\item maxdeg of any monomial in $\f_{\hat\gamma}(x)$ is $\leq$ (n-2).
\item If $\gamma$ is conjugate to a positive braid then 
 maxdeg($\f_{\hat\gamma}(x)$) is n-2 and all coefficients are positive
\item $\f_{\hat\gamma}(1)$ = w($\gamma$) 
(where w($\gamma$) is the writhe of $\gamma$).
\end{enumerate}
\end{teo}

\begin{proof} The proof for the case of braided knots is easier than the general case described by Fiedler.
One has to see that $\hat{\beta}\longmapsto\f_{\widehat{\beta}}(x)$ is a well defined map. 
More precisely I need to show firstly that
\[
\begin{array}{c}
\f_{\widehat{\alpha \beta}}(x)\eq \f_{\widehat{\alpha \sigma_j \sigma_j^{-1}\beta}}(x)\\
\f_{\widehat{\alpha \sigma_j \sigma_{j+1} \sigma_j \beta}}(x)
\eq \f_{\widehat{\alpha \sigma_{j+1} \sigma_j \sigma_{j+1} \beta}}(x)\\
\f_{\widehat{\alpha \sigma_j \sigma_i \beta}}(x)
\eq \f_{\widehat{\alpha \sigma_i \sigma_j \beta}}(x){\rm\ with\ }|i-j|\geq 2\\
\end{array}
\]

\noindent The last two equalities have to be proved because the relations in $\braid_n$ are 
$\sigma_j \sigma_{j+1} \sigma_j\eq \sigma_{j+1} \sigma_j \sigma_{j+1}$ for $j\eq\overline{1,n-2}$ and
$\sigma_j \sigma_i\eq \sigma_i \sigma_j $ for $|i-j|\geq 2$.
After smoothing the r-th crossing of $\widehat\beta$,
$\pi_\beta$ will break into a product of two disjoint cycles. 
Then m(r) is the length of the disjoint cycle which contains $i_r$, the index of the r-th ascending string. 
I will denote by $\p_\theta(x)$ all the monomials in the Fiedler polynomial corresponding 
to the letters in $\theta$, and also $m(\sigma_{i_r})\eq m(i_r).$

\[
\f_{\widehat{\alpha \sigma_j\sigma_j^{-1}\beta}}(x)\eq
\p_\alpha(x) + x^{2m(\sigma_j)-n} - x^{2m(\sigma_j^{-1})-n} + \p_\beta(x) 
\]

\noindent I have to show that $m(\sigma_j) = m(\sigma_j^{-1})$. 
When I am smoothing any of these two crossings the
associated permutations are the same: $\pi_{\alpha}s_j\pi_{\beta}$. 
So the decompositions into disjoint cycles, when smoothing either 
$\sigma_j$ or $\sigma_j^{-1}$
are the same. For $\sigma_j$ the ascending string is the $i_j$ string and for
$\sigma_j^{-1}$ is also $i_j$ viewed at the top of $\sigma_j$. 
So $m(\sigma_j) = m(\sigma_j^{-1})$ and from
here I obtain that:

\[
\f_{\widehat{\alpha \sigma_j \sigma_j^{-1}\beta}}(x)\eq 
\p_\alpha(x) + \p_\beta(x)\eq 
\f_{\widehat{\alpha\beta}}(x)
\]

\noindent In the same way one can prove the last two equalities, as well as the invariance
under conjugation. 
Everything else in the theorem can be proved using 
a "skein relation". The assertion is that the difference of the
respective values of the Fiedler polynomial on 
the braids $\alpha\sigma_j\beta$ and 
$\alpha\sigma_j^{-1}\beta$ is a symmetric polynomial 
( see \eqref{eq:fiedler} ). 
One observation is that the ascending string when smoothing
$\sigma_j^{-1}$ is the descending string for the case when 
smoothing $\sigma_j$, and so $m(\sigma_j^{-1})\eq n-m(\sigma_j)$.

\begin{equation}
\label{eq:fiedler}
\begin{split}
&\hskip 1in \f_{\widehat{\alpha \sigma_j \beta}}(x) - 
\f_{\widehat{\alpha \sigma_j^{-1}\beta}}(x)\eq\\
&\p_\alpha(x) + 
x^{2m(\sigma_j)-n} + \p_\beta(x) - [\p_\alpha(x) - 
x^{2m(\sigma_j^{-1})-n} + \p_\beta(x)]\eq\\
&\hskip 1in x^{2m(\sigma_j)-n} + x^{n-2m(\sigma_j)}\\
\end{split}
\end{equation}
The rest of the proof is exactly as in Fiedler's paper.
\end{proof}
\vskip .1in

\begin{rem}
From now on we can think of $\beta$ as being the conjugacy class of $\beta\in\braid_n$. 
The closure of any braid in $\beta$ is the same braided knot $\widehat\beta$. 
\end{rem}

Using only the Definition {\bf\ref{dfn:fiedler}} one can prove the following:

\begin{prop} 
\label{prop:fiedlerme}
The Fiedler polynomial is an order 1 invariant for braided
knots in $\er^3\setminus\{z-axis\}$.
\end{prop}

\begin{proof} This observation is implicit in Morton's work 
~\cite{MR2000i:57025}, using the well known theorem of 
Birman and Lin ~\cite{MR94d:57010}. A direct proof can be given by a direct calculation.

I use the symbol $S_i$ for the "singular" generator 
of singular braid monoid (see ~\cite{MR94b:57007}). 
That is an elementary braid in which 
the strings $i$ and $i+1$ intersect transversally at a point. 
Now I need to see that if I consider
a singular closed braid with two singularities $S_i$ and $S_j$ 
then the Fiedler polynomial vanishes. Solving each singularity in the two possible ways one gets:

\[
\f_{\widehat{A S_i B S_j C}}(x)\eq
\f_{\widehat{A \sigma_i B \sigma_j C}}(x) - \f_{\widehat{A \sigma_i B \sigma_j^{-1} C}}(x) -
\f_{\widehat{A \sigma_i^{-1} B \sigma_j C}}(x) + 
\f_{\widehat{A \sigma_i^{-1} B \sigma_j^{-1} C}}(x)].
\]

Now using the "skein relation" mentioned in the proof of the 
Theorem {\bf\ref{teo:fiedler}} I will obtain:

\[
\begin{array}{c}
\p_A(x) + x^{2m(i)-n} + \p_B(x) + x^{2m(j)-n} + \p_C(x)\\
\\
-[\p_A(x) + x^{2m(i)-n} + \p_B(x) - x^{n-2m(j)} + \p_C(x)]\\
\\
-[\p_A(x) - x^{n-2m(i)} + \p_B(x) + x^{2m(j)-n} + \p_C(x)]\\
\\
+ \p_A(x) - x^{n-2m(i)} + \p_B(x) - x^{n-2m(j)} + \p_C(x) = 0
\end{array}
\]
\end{proof}
\vskip .1in

\section{Exchange moves and Fiedler polynomial}
\label{section:excfiedler}

In this section all braids considered are representatives of knots.
I will apply Fiedler's polynomial to exchange related closed braids, 
$\widehat{\beta_1}$ and $\widehat{\beta_2}$. 
Recall from Definition {\bf \ref{def:exchange}}, that such braids 
look like this: $\beta_1\eq X\sigma_n^{-1}Y\sigma_n$ 
and $\beta_2\eq X\sigma_nY\sigma_n^{-1}$.

I want to discuss in some more details the braids $X$ and $Y$
which form the two exchange related braids.
If $X$ is a word only in letters $\sigma_1,\cdots,\sigma_{n-2}$ 
then $X$ commutes with $\sigma_n^{\pm1}$ 
and $\beta_1\eq\sigma_n^{-1}XY\sigma_n$, 
$\beta_2\eq\sigma_nXY\sigma_n^{-1}$ which are conjugate braids. 
Similarly if $Y$ is a word in letters 
$\sigma_1,\cdots,\sigma_{n-2}$
then $Y$ commutes with $\sigma_n^{\pm1}$ and both 
$\beta_i$'s are conjugate to $XY$. 
I am not interested in either of these two cases. 
So I will assume that $X,Y$ contain the letter $\sigma_{n-1}^{\pm1}$ in their expression. 

Suppose now $X = X_1\sigma_{n-1}^{2k}X_2$ with no 
other appearance of $\sigma_{n-1}$, which means that $X_1, X_2$
are words in letters $\sigma_1,\ldots ,\sigma_{n-2}$.
Then $\pi_X = \pi_{X_1}\pi_{X_2}$ and both 
$\pi_{X_1}$ and $\pi_{X_2}$ fix the letter $n$.
$\pi_{\beta_2}$ =  
$\pi_{X_1}\pi_{X_2}s_n\pi_Ys_n$ = $s_n\pi_{X_1}\pi_{X_2}\pi_Ys_n$ 
(because the multiplication of disjoint cycles is commutative). 

\noindent The last expression is conjugate to $\pi_{X_1}\pi_{X_2}\pi_Y$ 
which cannot be an $(n+1)$-cycle because it 
fixes the letter $n+1$. 
Also $X = X_1\sigma_{n-1}^{2k_1}X_2\sigma_{n-1}^{2k_2}\cdots
\sigma_{n-1}^{2k_l}X_l$, with $X_i\in\braid_{n-1}$,
with no other appearances of $\sigma_{n-1}$ 
wouldn't give us by closure a knot, because the n-string of the braid
would close itself to a component of the link. The same is true for $Y$. 
So at least one of the appearances of $\sigma_{n-1}$ in $X$ or $Y$, should be at an odd power. 

Consider also the braid $\beta = X\sigma_n Y\sigma_n$. The braids  
$\beta, \beta_1, \beta_2$ are simultaneously knots, that means their 
associated permutations are 
$\pi_\beta\eq\pi_{\beta_1}\eq\pi_{\beta_2}\eq$
$\pi_Xs_n\pi_Ys_n$ 
are $(n+1)$--cycles and I have the
relation w($X$) + w($Y$) $\equiv n$ (mod 2).
I want to see when $\f_{\widehat\beta_1}(x)$
 = $\f_{\widehat\beta_2}(x)$, for $\beta_1, \beta_2$ exchange related.
I will denote with $m_1(\sigma_n)$ the winding number of the ascending string
of the first $\sigma_n$ in the expression of $\beta$, and with $m_2(\sigma_n)$
the winding number of the ascending string for the second one.

\[
\begin{array}{c}
\f_{\hat\beta}(x) - \f_{\hat{\beta_1}}(x) 
\eq x^{2m_1(\sigma_n)-(n+1)} + x^{(n+1)-2m_1(\sigma_n)},\\
\f_{\hat\beta}(x) - \f_{\hat\beta_2}(x) 
\eq x^{2m_2(\sigma_n)-(n+1)} + x^{(n+1)-2m_2(\sigma_n)}
\end{array}
\]

\begin{equation}
\label{eq:skein}
\begin{array}{c}
\f_{\widehat\beta_1}(x) - \f_{\widehat\beta_2}(x) \eq \\ \\
x^{2m_2(\sigma_n)-(n+1)} + x^{(n+1)-2m_2(\sigma_n)} -
x^{2m_1(\sigma_n)-(n+1)} - x^{(n+1)-2m_1(\sigma_n)}
\end{array}
\end{equation}

Because $\pi_{\beta}\eq\pi_Xs_n\pi_Ys_n$ is a ($n+1$)--cycle
whenever I set $s_n\eq 1$ one of its appearances,
(smooth one of the $\sigma_n$ in the braid) I will get the 
resulting permutation to be a product of 2 disjoint cycles. 
So one of the cycles of $\pi_X\pi_Ys_n$ has length $m_1(\sigma_n)$,
and one of the cycles of $\pi_Xs_n \pi_Y$ has length $m_2(\sigma_n)$.

The above difference is a symmetric polynomial and depends only 
of the permutations $\pi_X\pi_Ys_n$ and $\pi_Xs_n \pi_Y$. Let $A\eq\pi_X$ and $B\eq s_n\pi_Ys_n$ then 
$\pi_X\pi_Ys_n = As_nB, \pi_Xs_n\pi_Y = ABs_n$.

\begin{rem}I will restrict my attention to the case when the above permutations $A$ and $B$ are both cycles, 
and I will study the case when they multiply to give us a full cycle which is $\pi_{\beta}$ the associated permutation
a braided knot $\hat{\beta}$.
\end{rem}
\vskip .1in

\begin{dfn}Let $l$ be the length of one of 
the cycles of $\pi_X\pi_Ys_n$. 
\end{dfn}

\noindent So is either $l$ is $m_1(\sigma_n)$ or $(n+1) - m_1(\sigma_n)$.

\begin{prop}
\label{prop:exrelength}
Let $n\geq 4$. If n+1 is odd and l$\eq\frac{n+2}{2}$ then $\f_{\widehat\beta_1}\eq\f_{\widehat\beta_2}$ for all
pairs of braids $\beta_1$ and $\beta_2$ related by an exchange move;
conversely, if  $\f_{\widehat\beta_1}\eq\f_{\widehat\beta_2}$ for all
pairs of braids $\beta_1$ and $\beta_2$ related by an exchange move then n+1 and l$\eq\frac{n+2}{2}$.
\end{prop} 

I mentioned the condition $n\geq 4$, (so the exchange related
braids will live in $\braid_5$ at least), because for braid index 3,
the exchange move does not change the conjugacy class of the braids 
and from the conclusion of the proposition we see that the values of
the Fiedler polynomial are different for two exchange related
braids living in $\braid_4$. 

To prove the above result I will study how the permutations {\it A} and $B$ multiplied give a full cycle. 
Let $A\eq (i_1 \cdots i_a)$, $B\eq (j_1,\cdots,j_b)$ be both cycles.
Consider $A\cup B\eq \{1, 2, \cdots, n+1\}$, and $A\cap B\eq \{t_1,\cdots,t_m\}$ to be the 
union and respectively the intersection of the sets
of digits appearing in $A$ and $B$.
Consider $t_1,\cdots,t_m$ to be written in an increasing order.
In both $A$ and $B$ these digits will appear in different order. Denote by $\nu_A$ the permutation
obtained by considering the digits $t_1,\cdots,t_m$ in the order they appear in $A$, and similarly $\nu_B$.
Denote also by $\nu_{A, B}\eq\nu_A^{-1}\cdot\nu_B$.
Let's see what these permutations look like calculating
a few examples (computations are done with ~\cite{GAP4}).

\flushleft{ 1. (1 7 3 9 5 4 2 8 6)(1 2)\eq (1 7 3 9 5 4)(2 8 6),
$A\cap B\eq (1\ 2),\ \nu_A\eq\nu_B\eq id,\ \nu_{A, B}\eq id$}

\flushleft{ 2. \quad (1 7 3 9 5 4 2 8 6)(1 2 3)\eq (1 7)(2 8 6)(3 9 5 4)
$A\cap B\eq (1\ 2\ 3),\ \nu_A\eq (2\ 3),\ \nu_B\eq id,\ \nu_{A, B}\eq (2\ 3)$}

\flushleft{ 3. \quad (1 7 3 9 5 4 2 8 6)(1 2 3 5 6)\eq (1 7 5 4 3 9 6 2 8)
$A\cap B\eq (1\ 2\ 3\ 5\ 6),\ \nu_A\eq (2\ 3\ 5),\ \nu_B\eq id,\ \nu_{A, B}\eq(2\ 5\ 3)$}

\flushleft{ 4. \quad (1 3 4 7 10 5 8 11)(6 7 5 3 2 9 8)\eq (1 2 9 8 11)(3 4 5 6 7 10)
$A\cap B\eq (3\ 5\ 7\ 8),\ \nu_A\eq (5\ 7),\ \nu_B\eq (3\ 7),\ \nu_{A, B}\eq (3\ 7\ 5)$}

\flushleft{ 5. (1 3 4 7 10 5 8 11)(6 7 5 3 2 9)\eq (1 2 9 6 7 10 3 4 5 8 11)
$A\cap B\eq (3\ 5\ 7),\ \nu_A\eq (5\ 7),\ \nu_B\eq (3\ 7),\ \nu_{A, B}\eq (3\ 7\ 5)$}

\flushleft{ 6. \quad (1 3 4 7 10 5 8 11)(6 7 2 3 9 5)\eq (1 9 5 8 11)(2 3 4)(6 7 10)
$A\cap B\eq (3\ 5\ 7),\ \nu_A\eq (5\ 7),\ \nu_B\eq (3\ 7\ 5),\ \nu_{A, B}\eq (3\ 7)$}

\noindent These computations lead to the following:

\begin{lema}
\label{lema:fullcycle}
Let A, B cycles, as above. Then $A\cdot B$ is a full cycle if and only if 
$|A\cap B|$ is odd and $\nu_{A, B}$ has an even number of inversions.
\end{lema}

\begin{proof}For the direct implication suppose first $A\cap B\eq \{t_1,\cdots,t_{2m}\}$ 
and $\nu_{A, B}$ is identity. We have 
$i_{k_p}\eq t_p\eq j_{k_p}$.

\[
\begin{array}{l}
A\cdot B\eq (i_1,\cdots,i_{k_1},\cdots,i_{k_{2m}},\cdots,i_a)\cdot
(j_1,\cdots,j_{k_1},\cdots,j_{k_{2m}},\cdots,j_b)\\
\hskip .32in\eq (i_1,\cdots,i_{k_1-1},j_{k_1+1}\cdots,j_{k_2},i_{k_2+1},
\cdots,j_{k_{2m}},i_{k_{2m}+1},\cdots,i_a)\cdot \\
\hskip .52in (j_1,\cdots,j_{k_1-1},j_{k_1},i_{k_1+1}\cdots,i_{k_2-1},j_{k_2+1},
\cdots,i_{k_{2m}-1},j_{k_{2m}+1},\cdots,j_b) 
\end{array}
\]

\noindent If we are changing $\nu_{A, B}$ with an even permutation, first we 
will encounter $i_a$ and this will give us a product of 2 cycles as above which is a contradiction.

Now suppose $A\cap B\eq \{t_1,\cdots,t_{2m+1}\}$ 
and sign($\nu_{A, B})\equiv$ 1 (mod 2).
I will consider the case $\nu_{A, B}$ = (1 2). So in this case we have 
$i_{k_1}\eq t_1\eq j_{k_2}$ and $i_{k_2}\eq t_2\eq j_{k_1}$ and for the 
rest of $p$'s $i_{k_p}\eq t_p\eq j_{k_p}$.

\[
\begin{array}{l}
A\cdot B\eq
(i_1,\cdots,i_{k_1},\cdots,i_{k_{2m+1}},\cdots,i_a)\cdot
(j_1,\cdots,j_{k_1},\cdots,j_{k_{2m+1}},\cdots,j_b)\\ 
\hskip .31in \eq (i_1,\cdots,i_{k_1-1},j_{k_2+1}\cdots,j_{k_3},i_{k_3+1},
\cdots,j_{k_{2m}},i_{k_{2m}+1},\cdots,i_a)\cdot \\
(i_{k_1},\cdots,i_{k_2-1},j_{k_1+1}\cdots,j_{k_2-1})\cdot
(i_{k_2},\cdots,i_{k_3-1},j_{k_3+1},\cdots,i_{k_{2m+1}-1},j_{k_{2m}+1},\cdots,j_b,j_1,\cdots,j_{k_1-1}) .
\end{array}
\]

Again if we are changing $\nu_{A,B}$ with an even permutation we will 
get the same kind of decomposition in disjoint cycles for $A\cdot B$ as above.

For the other implication I will consider 
$A\cap B\eq \{t_1,\cdots,t_{2m+1}\}$ and $\nu_{A,B}$ is identity.

\[
\begin{array}{l}
A\cdot B\eq
(i_1,\cdots,i_{k_1},\cdots,i_{k_{2m+1}},\cdots,i_a)\cdot
(j_1,\cdots,j_{k_1},\cdots,j_{k_{2m+1}},\cdots,j_b)\\
\hskip .32in \eq (i_1,\cdots,i_{k_1-1},j_{k_1+1}\cdots,j_{k_2},i_{k_2+1},
\cdots,i_{k_{2m+1}-1},j_{k_{2m+1}+1},\cdots,j_b,\\
\hskip .61in j_1,\cdots,j_{k_1},i_{k_1+1},\cdots,i_{k_2-1},j_{k_2+1},\cdots,
i_{k_{2m}-1},j_{k_{2m}+1},\cdots,i_a)\\
\end{array}
\]

Changing $\nu_{A,B}$ with an even permutation we 
get again an ($n+1$)-cycle because in the multiplication we 
are getting first $j_b$ and then at the end $i_a$.
\end{proof}
\vskip.1in

\noindent With $A, B$ as above, $A\cap B = \{t_1,\cdots,t_{2m+1}\}$, and $n,n+1$ are not in this set.
Then $As_nB\cap ABs_n = \{t_1,\cdots,t_{2m+1},n\}$ so the cardinal of this set is even.

\begin{lema}
\label{lema:permlength}
Let A, B as above, such that $A\cdot B$ is a full cycle,
then $\pi_X\pi_Ys_n$ and $\pi_Xs_n\pi_Y$ have 
disjoint cycle decomposition of lengths l,
n+1-l and respectively l-1, n+2-l. 

\end{lema}
\begin{proof} 
I will consider as in the previous proof only the permutation $\nu$
of the intersection set $\{t_1,\cdots,t_{2m-1}\}$ is the identity 
and $n$ and $n+1$ are at the end between $i_{2m-1}$ 
and $i_a$ and respectively between $j_{2m-1}$ and $j_b$. 
By cyclic permutations of digits we can arrange to have this
order in $A$. So the above assumption is only for $B$.

\begin{flushleft}
$A\eq (i_1,\cdots,i_{k_1},\cdots,i_{k_{2m-1}},\cdots,i_u,\cdots,i_a)$\quad
 where\quad 
$i_u\eq n$.
\smallskip
$B\eq (j_1,\cdots,j_{k_1},\cdots,j_{k_{2m-1}},\cdots,j_v,\cdots,j_b)$
\quad where \quad $j_v\eq n+1$.
\smallskip
$s_nB\eq (j_1,\cdots,j_{k_1},\cdots,j_{k_{2m-1}},
\cdots,j_{v-1},j_v,n,\cdots,j_b)$
\smallskip
$Bs_n\eq (j_1,\cdots,j_{k_1},\cdots,j_{k_{2m-1}},\cdots,j_{v-1},n,j_v,\cdots,j_b)$
\end{flushleft}

\[
\begin{array}{l}
As_nB\eq
(i_1,\cdots,i_{k_1},\cdots,i_{k_{2m-1}},\cdots,i_u,\cdots,i_a)\cdot
(j_1,\cdots,j_{k_1},\cdots,j_{k_{2m-1}},\cdots,j_{v-1},j_v,n,\cdots,j_b) \\

\hskip 0.36in\eq (i_1,\cdots,i_{k_1-1},j_{k_1+1}\cdots,j_{k_2},i_{k_2+1},\cdots,i_{k_{2m-1}-1},j_{k_{2m-1}+1},\cdots,j_{v-1},j_v,n,\cdots,i_a)
\cdot \\
(i_{k_1},\cdots,i_{k_2-1},j_{k_2+1}\cdots,j_{k_3},i_{k_3+1},\cdots,j_{k_{2m-1}},i_{k_{2m-1}+1},
\cdots,i_{u-1},j_{v+1},\cdots,j_b,j_1,\cdots,i_{k_1-1}) 
\end{array}
\]

\[
\begin{array}{l}
ABs_n\eq
(i_1,\cdots,i_{k_1},\cdots,i_{k_{2m-1}},\cdots,i_u,\cdots,i_a)\cdot 
(j_1,\cdots,j_{k_1},\cdots,j_{k_{2m-1}},\cdots,j_{v-1},n,j_v,\cdots,j_b)\\
\hskip .355in\eq (i_1,\cdots,i_{k_1-1},j_{k_1+1}\cdots,j_{k_2},i_{k_2+1},
\cdots,i_{k_{2m-1}-1},j_{k_{2m-1}+1},\cdots,j_{v-1},n,i_{u+1},\cdots,i_a)\cdot\\
(i_{k_1},\cdots,i_{k_2-1},j_{k_2+1}\cdots,
j_{k_3},i_{k_3+1},\cdots,j_{k_{2m-1}}, i_{k_{2m-1}+1},\cdots,i_{u-1},j_v,j_{v+1},\cdots,j_b,j_1,\cdots,i_{k_1-1})\\ 
\end{array}
\]

\vskip .1in
Sometimes the length of the first cycle in the decomposition 
of $As_nB$ will be $l-1$
and the length of the similar cycle in $ABs_n$ will be $l$. 
Such an example may be 
obtained as follows: let $i_u\eq i_a\eq n$ and $j_v\eq j_b\eq n+1$.
Consider also $A\eq (i_{k_1},i_{k_2},\cdots,i_{k_{2m-1}},\cdots,i_u),\ 
B\eq (j_1,j_{k_1},j_{k_2},\cdots,j_{k_{2m-1}},\cdots,j_{v-1},j_v).$

Then $s_nB\eq (j_{k_1},j_{k_2},\cdots,j_{k_{2m-1}},\cdots,j_{v-1},j_v,n)$,
and

$Bs_n\eq (j_{k_1},j_{k_2},\cdots,j_{k_{2m-1}},\cdots,j_{v-1},n,j_v)$ 
where $i_{k_p}\eq j_{k_p}.$  
\[
\begin{array}{l}
As_nB\eq
(i_{k_1},i_{k_2},\cdots,i_{k_{2m-1}},\cdots,i_u)\cdot
(j_{k_1},j_{k_2},\cdots,j_{k_{2m-1}},\cdots,j_{v-1},j_v,n)\\
\hskip .36in\eq (i_{k_1},j_{k_3},j_{k_5},\cdots,j_{k_{2m-1}},
\cdots,i_{k_{2m-1}+1},\cdots,i_{u-1})\cdot  \\
\hskip .55in (i_{k_2},j_{k_4},j_{k_6},\cdots,j_{k_{2m-2}},
i_{k_{2m-2}+1},\cdots,i_{k_{2m-1}-1},j_{k_{2m-1}+1},\cdots,j_{v-1},j_v,n) .
\end{array}
\]

\[
\begin{array}{l}
ABs_n\eq 
(i_{k_1},i_{k_2},\cdots,i_{k_{2m-1}},\cdots,i_u)\cdot
(j_{k_1},j_{k_2},\cdots,j_{k_{2m-1}},\cdots,j_{v-1},n,j_v)\\

\hskip .36in\eq(i_{k_1},j_{k_3},j_{k_5},\cdots,j_{k_{2m-1}},
\cdots,i_{k_{2m-1}+1},\cdots,i_{u-1},j_v)\cdot  \\
\hskip .56in (i_{k_2},j_{k_4},j_{k_6},\cdots,j_{k_{2m-2}},
i_{k_{2m-2}+1},\cdots,i_{k_{2m-1}-1},j_{k_{2m-1}+1},\cdots,j_{v-1},n) .
\end{array}
\]

The same type of decomposition in disjoint cycles will be 
obtained for any appropriate
$A$ and $B$, meaning the lengths of the disjoint cycles in $As_nB$ 
and $ABs_n$ will be as in the above cases.
\end{proof}

\vskip .1in
{\it Proof} of Proposition {\bf\ref{prop:exrelength}}: Since $l$ is 
either $m_1(\sigma_n)$ or $(n+1)-m_1(\sigma_n)$ it follows 
from lemma {\bf\ref{lema:permlength}} that the difference of the two polynomials
will vanish if and only if  $l\eq n+2-l\Leftrightarrow 2l\eq n+2$.
So from here we get two things: n must be even and 
$l\eq\frac{n+2}{2}$.
\hfill {$\Box$}
\newpage

\begin{exemplu}
\label{ex:braid1}
An infinite family of pairs of braids in $\braid_5\ (n\eq 4)$  
which are not conjugate and the Fiedler polynomial cannot distinguish them.  
\end{exemplu}

\hskip .21in Let 
\[
\beta_{1,k}\eq\sigma_3\sigma_2\sigma_1\sigma_4^{-1}
\sigma_3\sigma_2^{2k+1}\sigma_1\sigma_4
\quad
{\rm and}
\quad
\beta_{2,k}\eq\sigma_3\sigma_2\sigma_1\sigma_4
\sigma_3\sigma_2^{2k+1}\sigma_1\sigma_4^{-1}.
\]

\hskip .21in They are related by exchange moves for each $k$, 
and using the conjugacy algorithm (see ~\cite{prep1}) we see that 
they are not conjugate. Their associated permutations are 
\[
\pi_{\beta_{1,k}}\eq
\pi_{\beta_{2,k}}\eq (1\ 3\ 4\ 2\ 5);\quad X\eq\sigma_3\sigma_2\sigma_1,\quad
Y\eq\sigma_3\sigma_2^{2k+1}\sigma_1.
\] 
$\pi_X\eq (1\ 2\ 3\ 4)$, $s_4\pi_Ys_4 = (1\ 2\ 3\ 5)$,
and from here we get
$\pi_X\pi_Ys_4\eq (1\ 3)\cdot (2\ 5\ 4)$
and $\pi_Xs_4\pi_Y\eq (1\ 3\ 5)\cdot (2\ 4)$.

In this example $m_1(\sigma_n)\eq 5 - m_2(\sigma_n)\eq 3$.
So $l\eq\frac{n+2}{2}\eq\frac{6}{2}\eq 3$,
and their Fiedler polynomial are equal.
Similar examples can be obtained by 
multiplying each of the braids with a pure braid.
\hfill {$\Box$}

\vskip .1in
\begin{exemplu}
\label{ex:braid2}
Two sets of braids in $\braid_5\ (n\eq 4)$ which can be distinguished 
using the Fiedler polynomial. 
\end{exemplu}

\hskip .21in Let 
\[
\beta_1\eq\sigma_3\sigma_2\sigma_1\sigma_4^{-1}
\sigma_2\sigma_1\sigma_3\sigma_4 \quad
{\rm and}\quad
\beta_2\eq\sigma_3\sigma_2\sigma_1\sigma_4
\sigma_2\sigma_1\sigma_3\sigma_4^{-1}.
\]

\hskip .21in They are related by exchange moves. Their associated permutations are
 
\[
\pi_{\beta_1}\eq
\pi_{\beta_2}\eq (1\ 5\ 3\ 4\ 2);\quad X\eq\sigma_3\sigma_2\sigma_1,
\quad Y\eq\sigma_2\sigma_1\sigma_3.
\] 

$\pi_X\eq (1\ 2\ 3\ 4)$, $s_4\pi_Ys_4 = (1\ 2\ 3\ 5)$,
and from here we get:
$\pi_X\pi_Ys_4\eq (1\ 2\ 3\ 4)\cdot (1\ 2\ 5\ 4\ 3)\eq (1\ 5\ 4\ 2)$
and $\pi_Xs_4\pi_Y\eq (1\ 2\ 3\ 4)\cdot (1\ 2\ 4\ 5\ 3)\eq (1\ 4\ 2)\cdot (3\ 5).$

\noindent We have $m_1(\sigma_n)\eq 4; m_2(\sigma_n)\eq 2$.
In this case $l\eq 4$ or $l\eq 1\neq\frac{n+2}{2}\eq 3.$
Moreover the difference $\f_{\widehat\beta_1}(x) - \f_{\widehat\beta_2}(x)\eq 
x + x^{-1} - x^3 - x^{-3}\neq 0$. So we may conclude that they are not conjugate.
\hfill {$\Box$}

\vskip .1in
\begin{exemplu}
\label{ex:braid3} Let n be an even number. 
In this general case there are families of braids related 
by exchange moves which are distinguished by the 
Fiedler polynomial (if $l\eq\frac{n+2}{2}$).
\end{exemplu}
For $i\neq 1$, let
\[
\begin{array}{l}
\beta_1\eq\sigma_{n-1}\cdots\sigma_1\sigma_n^{-1}
\sigma_{n-1}\cdots\sigma_{i+1}\sigma_i\sigma_{i+1}^{-1}\cdots
\sigma_{n-1}^{-1}\sigma_n\quad {\rm and}\\ \\
\beta_2\eq\sigma_{n-1}\cdots\sigma_1\sigma_n
\sigma_{n-1}\cdots\sigma_{i+1}\sigma_i\sigma_{i+1}^{-1}\cdots
\sigma_{n-1}^{-1}\sigma_n^{-1}
\end{array}
\]
They are related by exchange moves for each $k$, 
Their associated permutations are: 
\begin{enumerate}
\item {for $i\eq 1\quad
\pi_{\beta_1}\eq
\pi_{\beta_2}\eq (1\ 2\cdots n-1\ n\ n+1)$ }  

\item {for $i\neq 1\quad
\pi_{\beta_1}\eq
\pi_{\beta_2}\eq (1\ 2\cdots i-1\ n+1\ i\cdots n)$ }
\end{enumerate}

\begin{center}
$X\eq\sigma_{n-1}\cdots\sigma_1$
and 
$\ Y\eq\sigma_{n-1}\cdots\sigma_{i+1}\sigma_i\sigma_{i+1}^{-1}\cdots
\sigma_{n-1}^{-1}$
\end{center}
$\pi_X\eq (1\ 2\cdots n)$, $s_n\pi_Ys_n = (i\ n+1)$,
and from here we get:

\begin{enumerate}
\item {
$i\eq 1$
\[
\begin{array}{l}
\pi_X\pi_Ys_n\eq (1\ 2\cdots n)\cdot (1\ n+1\ n) 
\eq (1\ 2\cdots n-1)\cdot (n\ n+1),\\
\pi_Xs_n\pi_Y\eq(1\ 2\cdots n)\cdot (1\ n\ n+1)
\eq (1\ 2\cdots n-1\ n+1)
\end{array}
\]}
and 
\item{ 
$i\neq 1$
\[
\begin{array}{l}
\pi_Xs_n\pi_Y\eq (1\ 2\cdots n)\cdot (i\ n+1\ n)\\
\hskip .56in \eq (1\ 2\cdots i-1\ n+1\ n)\cdot (i\ i+1\cdots n-1),\\
\pi_X\pi_Ys_n\eq (1\ 2\cdots n)\cdot (i\ n\ n+1)\\
\hskip .56in \eq (1\ 2\cdots i-1\ n)\cdot (i\ i+1\cdots n-1\ n+1)
\end{array}
\]
}
\end{enumerate}

Independent of $i$ we have that $m_1(\sigma_n)\eq i + 1$
and $m_2(\sigma_n)\eq n + 1 - i$,
so $l\eq i+1$, and the above permutations are 
conjugate if and only if $\frac{n+2}{2} = i+1 \Leftrightarrow
i\eq\frac{n}{2}$.
\hfill {$\Box$}
\vskip .1in

\section{Finite type invariants coming from Kauffman bracket}
\label{section:newinv}

In his paper ~\cite{MR0899057:57006}, Kauffman describes a purely combinatorial way of obtaining the 
Jones polynomial. His main construction is a regular isotopy invariant for links in $\er^3$ 
which was called later the Kauffman bracket.  

I will do the same construction in the solid torus and the 
object I'll obtain is an invariant of conjugacy classes of braids.
Let  $\ar\eq\ce[a^{\pm 1},x]$, where $a$ and x are two variables. 
The $TL_n$ is the Temperley-Lieb algebra of index $n$ which is generated by $1, e_1,\cdots e_{n-1}$ and has the 
following defining relations:

\[
\begin{array}{rclc}
e_i^2 & \eq & d\cdot e_i &\\
e_ie_{i\pm 1}e_i & \eq & e_i & \\
e_ie_j & \eq & e_je_i & \text{if}\ |i-j| > 1 \\
\end{array}
\]

Above, the coefficient $d\eq -a^2-a^{-2}$ .
\vskip .2in

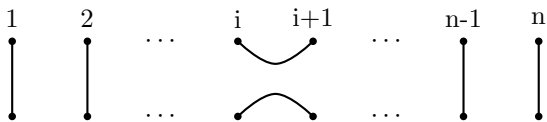
\begin{figure}[ht]
\begin{pspicture}(0,1.5)
\pscircle*(5,1){.05}
\pscircle*(5,0){.05}
\pscircle*(6,1){.05}
\pscircle*(6,0){.05}
\rput(7,1){$\cdots$}
\rput(7,0){$\cdots$}
\pscircle*(8,1){.05}
\pscircle*(8,0){.05}
\pscircle*(9,1){.05}
\pscircle*(9,0){.05}
\rput(10,1){$\cdots$}
\rput(10,0){$\cdots$}
\pscircle*(11,1){.05}
\pscircle*(11,0){.05}
\pscircle*(12,1){.05}
\pscircle*(12,0){.05}
\psline[linewidth=.8pt](5,1)(5,0)
\psline[linewidth=.8pt](6,1)(6,0)
\pscurve[linewidth=.8pt](9,1)(8.5,0.7)(8,1)
\pscurve[linewidth=.8pt](9,0)(8.5,0.3)(8,0)
\psline[linewidth=.8pt](11,1)(11,0)
\psline[linewidth=.8pt](12,1)(12,0)
\rput(5,1.3){1}
\rput(6,1.3){2}
\rput(8,1.3){i}
\rput(9,1.3){i+1}
\rput(11,1.3){n-1}
\rput(12,1.3){n}
\end{pspicture}
\caption{Generator $e_i$}
\label{fig:gener}
\end{figure}

\vskip .1in
\hskip .21in The elements of the Temperley-Lieb algebra can be viewed as
pictures. Consider $2n$ points lying on two parallel horizontal lines 
in the plane, $n$ on each line. The elements of the Temperley-Lieb algebra, which are also called states 
can be interpreted as a collection of $n$ disjoint arcs connecting these $2n$ points. 

Figure {\bf\ref{fig:gener}} shows a picture of the 
algebra generator $e_i$. These generators can be multiplied by putting one in top of the other to give all
the elements of the algebra. The identity is the set of segments connecting the $i$-th point on the top line with
the $i$-th point on the bottom line.

\hskip .21in Consider the following map: 

\begin{equation}
\Phi:\braid_n \hookrightarrow\ce\braid_n\to  TL_n
\label{dia:inv}
\end{equation}

The map \eqref{dia:inv} is the representation 
of the braid group on the Temperley-Lieb algebra, given by: 

\[
\Phi(\sigma_i)\eq a^{-1}\cdot e_i + a\cdot 1, \quad \Phi(\sigma_i^{-1})\eq a\cdot e_i + a^{-1}\cdot 1.
\]
  
\hskip .21in Let me explain the "closure" operation.
Think the plane in which the states live as being the sectional plane of a cylinder as in the case of braids.
The closure of a state is the exactly the same as the one for braids, meaning connect $i$-th point on the top line
around the axis with the $i$-th point on the bottom line. 
Closing a state give us a set of unoriented curves in the solid 
torus ( the one winding around the axis considered), some curves being
contractible, some of them winding around the axis. 

\hskip .21in Let $p\eq\#$ of contractible components and $q\eq\#$ of homotopycally non-trivial components in the closure of a state.


\hskip .21in $TL_n$ is a vector space over $\ce$ with the basis given by the set of states. 
\begin{dfn} Let $f: TL_n\to\ar \quad f(w)\eq d^px^q $
for $w$ a state and extend by linearity to the entire vector space.
\end{dfn}

$f\circ\Phi$ associates to each braid a polynomial in $a^{\pm 1}$ and $x$. 
One can see that $f(e_i)\eq d\cdot x^{n-2}$ and $f(1)\eq x^n$.

\begin{lema}
\label{lema:trace}
$f: TL_n\to\ar$ is a trace function, in other
words, it has the following three properties for $(\forall)\ v,\ w\in TL_n$
and for $(\forall)\ \alpha\in\ce[a^{\pm 1}]$.
\begin{enumerate}
\item {$f(v + w)\eq f(v) + f(w)$}
\item {$f(\alpha\cdot v)\eq \alpha\cdot f(v)$}
\item {$f(v\cdot w)\eq f(w\cdot v)$} 
\end{enumerate}
\end{lema}

\begin{proof} The first two properties are true because we 
have defined the function $f$ as a linear extension. For last property 
let's look at two generators of the Temperley-Lieb algebra, $e_i, e_j$. 
Consider the closure of $e_ie_j$, we can isotope $e_j$ through
the closure arcs in the top of $e_i$. So the two closures are 
isotopic in the solid torus, so they have the same number of 
components, which means that $f(e_ie_j)\eq f(e_je_i)$.
From here, by the same argument we find that $f$ has
the third property for every element in the basis of the 
Temperley-Lieb algebra. Now since each word is a sum of states with 
some coefficients, by linearity we get the property in general.
\end{proof}
\vskip .1in

For father computations it will be useful to have some formulas 
for $\Phi(\sigma_i^k).$

\begin{lema}
\label{lema:image}
$\Phi(\sigma_i^k)\eq p_k(a)e_i+a^k$, where 
\[
p_k(a)\eq
\left\{
\begin{array}{lcr}
\sum_{l\eq 0}^{k-1} (-1)^l a^{k-2-4l}&if&k>0\\
1&if&k\eq 0\\
\sum_{l\eq 0}^{-k-1} (-1)^l a^{k+2+4l}&if&k<0\\
\end{array}
\right.
\]
\end{lema}

\begin{proof} Consider $k>0$. The proof will be by 
induction. For $k\eq 1$ we have $p_1(a)\eq a^{-1}$ which
is exactly what the formula give us. Assume the formula is
true for $k$. For $k+1$ we have to prove that  
\[
p_{k+1}(a)\eq\sum_{l\eq 0}^k (-1)^l a^{k+1-2-4l}
\]
\[
\Phi(\sigma_i^{k+1})\eq (a^{-1}e_i+a)\cdot (p_k(a)e_i+a^k)
\eq (a^{-1}p_k(a)d + a^{k-1} + ap_k(a))
\]
\[
\eq (-a^{-3}p_k(a) + a^{k-1})e_i + a^{k+1}
\]
The coefficient of $e_i$ is $p_{k+1}(a)$
\[
\eq (-a^{-3})\sum_{l\eq 0}^{k-1} (-1)^l a^{k-2-4l} + a^{k-1}
\eq \sum_{l\eq 0}^{k-1} (-1)^{l+1} a^{k-3-2-4l} + a^{k+1-2}
\]
\[
\eq \sum_{l\eq 0}^{k-1} (-1)^{l+1} a^{k+1-2-4(l+1)} + a^{k+1-2}
\]
\[
\eq \sum_{l\eq 1}^k (-1)^l a^{k+1-2-4l} + a^{k+1-2}
\eq \sum_{l\eq 0}^k (-1)^l a^{k+1-2-4l}.
\]
The case $k<0$ is similar. 
\end{proof}
\vskip .1in

\begin{prop}
\label{prop:invme}
Let $a\eq e^t$ in $\Phi(\beta)$. Then the coefficient of $t^k$
in $\Phi(\beta)$ is a finite order invariant of order k.
\end{prop}

\begin{proof} Consider a singular braid with $k+1$ transverse
double points, and consider its image in the Temperley-Lieb algebra,
through the canonical extension of $\Phi$ to singular braids.

I will show that all the coefficients of $t^k$
up to power k+1 obtained after we expand $a\eq e^t$, will be equal to zero.
Let $\beta\eq X_1\tau_{i_1}X_2\tau_{i_2}\cdots X_{k+1}\tau_{i_{k+1}}X_{k+2}$, 
where $X_i's$ are braids in $\braid_n$ and $\tau_i$ is the standard singular 
i-th braid generator (the strings i and i+1 cross transversally in
one point). 

\begin{equation}
\label{eq:Kfinite}
\begin{split}
 & \hskip 1in \Phi(X_1\tau_{i_1}X_2\tau_{i_2}\cdots X_{k+1}\tau_{i_{k+1}}X_{k+2})\eq \\
 & \\
 & \sum_{\epsilon_{i_1},\cdots,\epsilon_{i_{k+1}};{\epsilon_{i_l}\in\{\pm 1\}}}
\epsilon_{i_1}\ldots\epsilon_{i_{k+1}}
\Phi(X_1\sigma_{i_1}^{\epsilon_{i_1}}X_2\sigma_{i_2}^{\epsilon_{i_2}}
\cdots X_{k+1}\sigma_{i_{k+1}}^{\epsilon_{i_{k+1}}}X_{k+2})\eq \\
 & \eq\sum_{\epsilon_{i_1},\cdots,\epsilon_{i_{k+1}};{\epsilon_{i_l}\in\{\pm 1\}}}
\epsilon_{i_1}\ldots\epsilon_{i_{k+1}}
\Phi(X_1)(a^{-\epsilon_{i_1}}e_{i_1}+a^{\epsilon_{i_1}})
\cdots \Phi(X_{k+1})(a^{-\epsilon_{i_{k+1}}}e_{i_{k+1}}+a^{\epsilon_{i_{k+1}}})\Phi(X_{k+2})
\end{split}
\end{equation}

In the sum \eqref{eq:Kfinite} we have $2^{k+1}$ terms of the form
\[\epsilon_{i_1}\ldots\epsilon_{i_{k+1}}\Phi(X_1)(a^{-\epsilon_{i_1}}e_{i_1}+a^{\epsilon_{i_1}})
\cdots \Phi(X_{k+1})(a^{-\epsilon_{i_{k+1}}}e_{i_{k+1}}+a^{\epsilon_{i_{k+1}}})\Phi(X_{k+2})\] 
and expanding each parenthesis we will get for each of them a sum of 
$2^{k+1}$ terms, each one of them being like  
$\Phi(X_1)\cdot e_{i_1}\cdots \Phi(X_{k+1})\cdot e_{i_{k+1}}\cdot\Phi(X_{k+2})$, or with some $e_j$ missing.
I am interested in the coefficient of each such term. 
In fact all the coefficients of any of the terms are equal up to a sign.
Such term is realized $2^{k+1}$ times. The sign in front
of each appearance will be $(-1)^l$ where $l$ is the number of -1's in the $k+1$-tuple $\epsilon_{i_1},\ldots,\epsilon_{i_{k+1}}$.
For a given such $l$ we have $\binom{k+1}{l}$ ways of choosing $l$ -1's out of $k+1$ numbers. 
As for the $a$ factor for a given $l$ will be of the form $a^l\cdot (a^{-1})^{k+1-l}$, because 
$\Phi(\sigma_j^{-1})\eq (ae_i + a^{-1})$. 
Putting everything together we get that the coefficient in front of   
$\Phi(X_1)\cdot e_{i_1}\cdots \Phi(X_{k+1})\cdot e_{i_{k+1}}\cdot
\Phi(X_{k+2})$ is

\begin{equation}
\label{exp:1}
\sum_{l=0}^{k+1}
\binom{k+1}{l}
(-1)^la^l\cdot (a^{-1})^{k+1-l}\eq
(a^{-1}-a)^{k+1}
\end{equation}

For $a\eq e^t$, we have  $a^{-1}-a\eq (1-t+\cdots) - (1+t+\cdots)\eq -2t+\cdots$, 
and we will have $(-a+a^{-1})^{k+1}\eq (-2t)^{k+1}+\cdots\equiv_k 0$, 
so the coefficients of $t^j\eq 0$ for all $1\leq j\leq k$. Here $\equiv_k$ means truncating the terms of degree $\geq k+1$.
Since the coefficient of all states which appear in the decomposition 
\eqref{eq:Kfinite} is up to sign equal with \eqref{exp:1} 
we get the conclusion of the proposition.
\end{proof}
\vskip .1in

\begin{dfn}
\label{dfn:newinv}
Consider $f\circ\Phi(\beta)$, and let $a\eq e^t$.
The coefficient of $t^k$, a polynomial in x, is denoted by $\q_{\hat\beta,k} (x)$. 
\end{dfn}

\begin{corolar}
\label{corolar:newfiniteinv}
$\q_{\hat\beta,k} (x)$ is a $k-th$ order invariant.
\end{corolar}

\begin{proof} Consider the image of the word in
\eqref{eq:Kfinite} after replacing $a$ by $e^t$, through $f$, 
We have that $\q_{\hat\beta,k} (x)\eq 0$, because all the 
coefficients up to order $k+1$ are zero.
\end{proof}
\vskip .1in

\begin{prop}
\label{prop:invdifer}
$\q_{\hat\beta,1}(x)\neq\f_{\hat\beta}(x)$  
\end{prop}

\begin{proof} Consider $\sigma_1^3\in \braid_2$. Its closure is the 
right hand trefoil. $\f_{\widehat{\sigma_1^3}}(x)\eq 3$.
Consider $\Phi(\sigma_1^3)\eq (a^{-1}e_1+a)^3\eq p_3(a)e_1+a^3$ where 
$p_3(a)\eq a - a^{-3} + a^{-7}$. Then 
$f\circ\Phi(\sigma_1^3)\eq p_3(a)\cdot 
f\circ\Phi(e_1)+a^3\cdot f\circ\Phi(1)\eq p_3(a)\cdot d+a^3\cdot x^2.$
Let $a\eq e^t$. Looking only to the terms up to order 2 we will
get $d\equiv_2 -(1+2t)-(1-2t)\eq -2$, 
$p_3(a)\equiv_2 1+t - (1-3t) + (1-7t)\eq 1-3t$,
$a^3\eq 1+3t$.
As a result
$f\circ\Phi(\sigma_1^3)\equiv_2 (1-3t)\cdot(-2) + (1+3t)x^2\eq 
-2+x^2 + 3(2+x^2)t$, and $\q_{{\widehat{\sigma_1^3}},1}(x)\eq 3(2+x^2)$.
\end{proof}
\vskip .1in

\hskip .21in We see that the $\q_{\hat\beta,k}$ exists for any $ k$, and moreover 
can be defined for any braids, not only for braided knots.

\hskip .21in I will start the study of this new invariant on exchange related knots with a 
discussion on Morton's braid representative of the unknot appearing in \cite{MR84m:57006}.
A conjugate of it is  
\[
\beta\eq
\sigma_2\sigma_2\sigma_2\sigma_1^{-1}\sigma_2\sigma_3^{-1}
\sigma_2^{-1}\sigma_2^{-1}\sigma_1\sigma_2^{-1}\sigma_3
\]

\hskip .21in It belongs to $\braid_4$. To destabilize to the canonical braid representative of the unknot 
in $e\in\braid_1$, one needs either to stabilize first to $\braid_5$, 
as Morton shows in \cite{MR84m:57006}, or 
needs to use exchange moves as in \cite{MR92g:57010b}. 
In fact it is sufficient to use only one exchange move, not two as in \cite{MR92g:57010b}.
The notations for this example only, are: the standard generators of the braid
group, $\sigma_i$ are replaced with $i$, and $\sigma_i^{-1}$ with $\overline{i}$), a twiddle means that the
transformation is conjugation in $\braid_n$, and an arrow with no 
specification means that we are using only the braid relations, the rest being clear. 
\[
\begin{array}{ccccccc}
222\overline 1 2\overline {322}1\overline 2 3 & exchange & 22212\overline {32212}3 & 
\to & 22212\overline {32121}3 & \to & \\
22212\overline {31211}3 & \to & 22212\overline {132}3\overline{11} & 
\sim & \overline{11}22212\overline {1}2\overline {32} & \to & \\
\overline{211}22212\overline {1}2\overline {3} & destabilization & \overline{211}22212\overline {1}2 &
\sim & \overline{11}22212\overline {1} & \to & \\
\overline{11}22121\overline 1 & \to & \overline{11}2212 & 
\to & \overline{11}2121 & \sim & \\
\overline{1}212 & \to & \overline{1}121 & 
\sim & 12 & destabilization & e \\
\end{array}
\]

\hskip .21in Now let me consider the exchange related braids:

\[
\begin{array}{lr}
\beta_1\eq \sigma_2^{-1}\sigma_3\sigma_2\sigma_2\sigma_2
\sigma_1^{-1}\sigma_2\sigma_3^{-1}\sigma_2^{-1}\sigma_2^{-1}\sigma_1,
&  
\beta_2\eq \sigma_2^{-1}\sigma_3\sigma_2\sigma_2\sigma_2
\sigma_1\sigma_2\sigma_3^{-1}\sigma_2^{-1}\sigma_2^{-1}\sigma_1^{-1},
\end{array}
\] 
with 
$X\eq \sigma_2^{-1}\sigma_3\sigma_2\sigma_2\sigma_2$ and 
$Y\eq \sigma_2\sigma_3^{-1}\sigma_2^{-1}\sigma_2^{-1}.$ 
First let's calculate $\q_{\widehat\beta_1,k}(x)$, and I will start $\Phi(\beta_1)$. 

\begin{equation}
\label{eq:beta1}
\begin{split}
\Phi(\beta_1)\eq & \Phi(X)(ae_1+a^{-1})\Phi(Y)\Phi(a^{-1}e_1+a)\eq \\
&\eq\Phi(X)e_1\Phi(Y)e_1 + a^2\cdot\Phi(X)e_1\Phi(Y) + \\
&+a^{-2}\cdot\Phi(X)\Phi(Y)e_1 + \Phi(X)\Phi(Y). \\
\end{split}
\end{equation}

\hskip .21in The computations involve lemma {\bf\ref{lema:image}} and we will get:
\vskip .1in

$\Phi(X)e_1\Phi(Y)e_1\eq 
(-a^{-13}+2a^{-9}-3a^{-5}+2a^{-1}-a^3)\cdot e_3e_2e_1 + 
(-2a^{-1}+3a^3-2a^7+a^{11})\cdot e_2e_3e_1 + 
(-2a^{-3}+2a-a^5)\cdot e_3e_1 + 
(-a^{-11}+3a^{-7}-5a^{-3}+4a-2a^5)\cdot e_2e_1 + 
(a+2a^3-a^7)\cdot e_1 $
\vskip .1in

$\Phi(X)e_1\Phi(Y)\eq (a^{-5}-3a^{-1}+5a^3-3a^7+a^{11})\cdot e_2e_3e_1e_2 +
(a^{-7}-2a^{-3}+4a+a^3-a^5-a^7-a^9)\cdot e_1e_3e_2 + 
(a^{-11}-a^{-7}+a^{-3})\cdot e_3e_2e_1 + 
(a^{-7}-2a^{-3}+3a-a^5)\cdot e_2e_1e_3 + a\cdot e_1e_2e_3 + 
(a^{-9}-a^{-5}+2a^{-1})\cdot e_1e_3 + 
(a^{-9}-3a^{-5}+3a^{-1}-a^3)\cdot e_3e_2 +
(a^{-9}-2a^{-5}+2a^{-1})\cdot e_2e_3 + 
(a^{-9}-2a^{-5}+2a^{-1})\cdot e_2e_1 +
(2a^3-a^7)\cdot e_1e_2 + 
(a^{-11}-a^{-7}+a^{-3})\cdot e_3 +
(2a^{-7}-5a^{-3}+6a-2a^5)\cdot e_2 +
a\cdot e_1$ 
\vskip .1in

$\Phi(X)\Phi(Y)e_1\eq 
(-a^{-11}+a^{-7}-a^{-3}+2a+a^3-a^7-a^9)\cdot e_3e_2e_1 + 
(-a^{-11}+2a^{-7}-2a^{-3}+2a-a^5)\cdot e_2e_3e_1 +
(-a^{-13}+a^{-9}-a^{-5}+a^{-1})\cdot e_1e_3 + 
(-a^{-9}+2a^{-5}-2a^{-1}+3a^3-2a^7+a^{11})\cdot e_2e_1 +
a\cdot e_1 $ 
\vskip .1in

$\Phi(X)\Phi(Y)\eq 
(-a^{-11}+a^{-7}-a^{-3}+a)\cdot e_3e_2 + 
(-a^{-11}+2a^{-7}-2a^{-3}+2a-a^5)\cdot e_2e_3 +
(-a^{-13}+a^{-9}-a^{-5}+a^{-1})\cdot e_3 + 
(-a^{-9}+2a^{-5}-2a^{-1}+3a^3-2a^7+a^{11})\cdot e_2 + a$ 
\vskip .1in

\hskip .21in The image of \eqref{eq:beta1} through $f$. 
Let's start with the images of the states appearing:
\[
\begin{array}{c}
f(e_2e_3e_1e_2)\eq f(e_1e_3)\eq d^2\\
f(e_1e_3e_2)\eq f(e_3e_2e_1)\eq f(e_2e_1e_3)\eq f(e_1e_2e_3)\eq d\\
f(e_3e_2)\eq f(e_2e_3)\eq f(e_2e_1)\eq f(e_1e_2)\eq x^2\\
f(e_3)\eq f(e_2)\eq f(e_1)\eq d\cdot x^2
\end{array}
\]

\hskip .21in Replacing all these in \eqref{eq:beta1},we get:

\[
\begin{array}{c}
f\circ\Phi(\beta_1)\eq f(\Phi(X)e_1\Phi(Y)e_1) + 
a^2\cdot f(\Phi(X)e_1\Phi(Y)) + \\
a^{-2}\cdot f(\Phi(X)\Phi(Y)e_1) + f(\Phi(X)\Phi(Y))\eq \\ 
(-a^{-19}+2a^{-15}-2a^{-11}-a^{-7}+3a^{-3}-a^{-1}-5a-a^3 +\\
2a^5+a^7+4a^9+a^{11}-a^{13}+a^{17}) + \\
(a^{-15}-4a^{-11}+9a^{-7}-15a^{-3}-a^{-1}+9a-a^3-11a^5+4a^9-a^{13})\cdot x^2
+ a\cdot x^4
\end{array}
\]

and similarly

\[
\begin{array}{c}
f\circ\Phi(\beta_2)\eq f(\Phi(X)e_1\Phi(Y)e_1) + 
a^{-2}\cdot f(\Phi(X)e_1\Phi(Y)) + \\
a^2\cdot f(\Phi(X)\Phi(Y)e_1) + f(\Phi(X)\Phi(Y))\eq \\ 
(-a^{-3}-a^{-1}-a-a^3+a^5+a^7+2a^9+a^{11}+a^{13}) + \\
(-2a^{-3}-a^{-1}-4a-a^3-2a^5)\cdot x^2
+ a\cdot x^4
\end{array}
\]

\hskip .21in It is clear, even from the images in the Temperley-Lieb
algebra, that the two braids are different, but we can see
that also their images through $f$ are different. We can look at 
$\q_{\widehat{\beta_i},k}(x)$, where $i\eq 1,2$ and $k\geq 0.$
For example $\q_{\widehat{\beta_1},1}(x)\eq 70 - 14x^2 + x^4$ and
$\q_{\widehat{\beta_2},1}(x)\eq 54 - 10x^2 + x^4.$

\hskip .21in Consider now in general two braids which are exchange related 
$\beta_1\eq X\sigma_n^{-1}Y\sigma_n$ and 
$\beta_2\eq X\sigma_nY\sigma_n^{-1}.$
We would be interested to compute the difference 
$f\circ\Phi(\beta_1) - f\circ\Phi(\beta_2)$. Let's look first 
to the difference of the images of $\beta_i$ in the $TL_n$ algebra. 
\[
\begin{array}{c}
\delta\eq \Phi(X)\Phi(\sigma_n^{-1})\Phi(Y)\Phi(\sigma_n) -
\Phi(X)\Phi(\sigma_n)\Phi(Y)\Phi(\sigma_n^{-1})\eq \\
\Phi(X)\cdot(ae_n+a^{-1})\Phi(Y)\Phi(a^{-1}e_n+a) - 
\Phi(X)\cdot(a^{-1}e_n+a)\Phi(Y)\Phi(ae_n+a^{-1}) \\ 
\end{array}
\]

\hskip .21in So we can rewrite $\delta$ as:
\begin{equation}
\label{eqn:diff}
\delta\eq (a^2 - a^{-2})\cdot [\Phi(X)e_n\Phi(Y) - \Phi(X)\Phi(Y)e_n]
\end{equation}

\hskip .21in Using the power series expansion $a\eq e^t$, we see that 
$a^2 - a^{-2}\equiv_2 4t$. In case we want to look only
to $\q_{\widehat{\beta_1},1}(x) - \q_{\widehat{\beta_2},1}(x)$,
we have then to compute only the free term of 
$f(\Phi(X)e_n\Phi(Y)) - f(\Phi(X)\Phi(Y)e_n).$

\hskip .21in Let me investigate the braids in example \eqref{ex:braid1} using these invariants. 

\begin{exemplu}[Example\eqref{ex:braid1} revisited]
\label{ex:Kbraid1}
The following two sets of braids in $\braid_5$ are not conjugate, as we already know,
and they cannot be distinguished either by the Fiedler polynomial, or by 
$\q_{\widehat{\beta},1}(x)$, but they can be distinguished 
by $\q_{\widehat{\beta},2}(x)$
\[
\beta_{1,j}\eq\sigma_3\sigma_2\sigma_1\sigma_4^{-1}
\sigma_3\sigma_2^{2j+1}\sigma_1\sigma_4
\quad
{\rm and}
\quad
\beta_{2,j}\eq\sigma_3\sigma_2\sigma_1\sigma_4
\sigma_3\sigma_2^{2j+1}\sigma_1\sigma_4^{-1}.
\]
\end{exemplu}

\hskip .21in Set $X\eq\sigma_3\sigma_2\sigma_1$ and $Y\eq\sigma_3\sigma_2^{2j+1}\sigma_1$.
\[
\begin{array}{c}
\Phi(X)e_4\Phi(Y) - \Phi(X)\Phi(Y)e_4\eq \\
\Phi(\sigma_3)\Phi(\sigma_2\sigma_1)e_4 
\Phi(\sigma_2^{2j+1}\sigma_1) - 
\Phi(\sigma_3)\Phi(\sigma_2\sigma_1) 
\Phi(\sigma_2^{2j+1}\sigma_1)e_4\eq \\
a^{-2}\cdot e_3\Phi(\sigma_2\sigma_1)e_4e_3\Phi(\sigma_2^{2j+1}\sigma_1) -
a^{-2}\cdot e_3\Phi(\sigma_2\sigma_1)e_3\Phi(\sigma_2^{2j+1}\sigma_1)e_4 + \\
\Phi(\sigma_2\sigma_1)e_4e_3\Phi(\sigma_2^{2j+1}\sigma_1) -
\Phi(\sigma_2\sigma_1)e_3\Phi(\sigma_2^{2j+1}\sigma_1)e_4
\end{array}
\]

\hskip .21in The last equality is obtained using the expression for $\Phi(\sigma_3)$ and
the commutativity relations of $e_4$ with $e_1, e_2$.
We can reduce it more because we are interested in the images through
$f$ which are the same for cyclic permutations of
the words. So in the image of $f$ the last
difference will cancel out. This lead us to the investigation of:
\[
\begin{array}{c}
\label{eq:diff10}
a^{-2}\cdot \Phi(\sigma_2\sigma_1)e_3\Phi(\sigma_2^{2j+1}\sigma_1)
\cdot [e_3e_4 - e_4e_3]\eq \\

\eq 
a^{-1}p_{2j+1}(a)\cdot (e_2e_1e_3e_4 - e_2e_1e_3e_2e_4e_3) + \\
\left(a^{-3}p_{2j+1}(a) + a^{2j-2}d + 2a^{2j}\right)\cdot
(de_2e_1e_3e_4 - e_2e_1e_3) + \\
ap_{2j+1}(a)\cdot (e_1e_3e_4 - e_1e_3e_2e_4e_3) + 
ap_{2j+1}(a)\cdot (e_3e_1e_4 - e_3e_2e_4e_3e_1) + \\
\left(a^{-1}p_{2j+1}(a) + a^{2j}d + 2a^{2j+2}\right)\cdot 
(de_1e_3e_4 - e_1e_3) + \\ 
a^{2j+2}\cdot (de_2e_3e_4 - e_2e_3) + \\
a^{-1}p_{2j+1}(a)\cdot (e_2e_1e_3e_4 - e_2e_1e_4e_3) + \\
a^{2j+4}\cdot (de_3e_4 - e_3) + 
ap_{2j+1}(a)\cdot (e_2e_3e_4 - e_2e_4e_3) 
\end{array}
\]

\hskip .21in Let's evaluate the function $f$ for the above states.
\[
\begin{array}{c}
f(e_2e_1e_3e_4)\eq f(e_2e_1e_3e_2e_4e_3)\eq x \\
f(e_1e_3e_2e_4e_3)\eq f(e_3e_2e_4e_3e_1)\eq dx \\
f(e_2e_1e_3)\eq f(e_3e_1e_4)\eq 
f(e_2e_3e_4)\eq f(e_2e_4e_3)\eq dx \\
f(e_1e_3)\eq f(e_3e_2e_4e_3)\eq d^2x;
f(e_3e_4)\eq f(e_2e_3)\eq x^3; f(e_3)\eq dx^3 \\
f(e_2e_1e_3e_4)\eq f(e_2e_1e_4e_3)\eq x
\end{array}
\]

\hskip .21in The difference $f(\Phi(X)e_4\Phi(Y)) - f(\Phi(X)\Phi(Y)e_4)$
using the expressions above becomes:
\[
a^3p_{2j+1}(a)\cdot (x^3 - d^2x) +
a^{2j+2}\cdot (d^2x - x^3),
\]

\hskip .21in The image of $\delta$ through $f$, will be: 

\[
\begin{array}{c}
f(\delta)\eq (a^2 - a^{-2})\cdot
[f(\Phi(X)e_4\Phi(Y)) - f(\Phi(X)\Phi(Y)e_4)]\eq \\
(a^2 - a^{-2})\cdot [(a^3p_{2j+1}(a) - a^{2j+2})\cdot x^3 +
d^2(a^{2j+2} - a^3p_{2j+1}(a))\cdot x]
\end{array} 
\]

\hskip .21in and after expanding $a\eq e^t$, we get up to degree 3 in $t$:

\[
f(\delta)\equiv_3 4\left(-4jx^3 + 16jx\right)\cdot t^2.
\]

\hskip .21in So we have learned that  
$\q_{\widehat{\beta_1,j},k}(x)-\q_{\widehat{\beta_2,j},k}(x) = 0$
if $k\eq 0, 1$, and the first non-zero difference is
$\q_{\widehat{\beta_1,j},2}(x)-\q_{\widehat{\beta_2,j},2}(x) = -16j(x^3 - 4x).$
\hfill{$\Box$}

Because of the above example and similar computations I will make the following:
\begin{conj}
\label{conj:new}
$\q_{\widehat\beta,1}(x)$ braids will vanish in the same way
on exchange related as the Fiedler's polynomial does
(see Proposition {\bf\ref{prop:exrelength}}).  
\end{conj}

{\it E-mail address:} {\tt radu.popescu@imar.ro}

\end{document}